\documentclass[12pt]{article}
\usepackage{amsxtra,amssymb,amsthm,amsmath,latexsym}

\theoremstyle{plain}
\newtheorem{theorem}{Theorem}[section]
\numberwithin{equation}{section}

\def\R{{\mathbb R}}

\def\bee{\begin{equation*}}
\def\eee{\end{equation*}}

\def\oH{\buildrel\circ\over H}
\def\oH1{\buildrel\circ\over H\kern-.02in{}^1}

\def\loc{{\hbox{\,loc\,}}}

\def\be{\begin{equation}}
\def\ee{\end{equation}}

\begin{document}

\newcommand{\coloneqq}{\mathrel{\raise.095ex\hbox{:}\mkern-4.2mu=}}
\def\ecolon{\mathrel{\mathop:}=}

\newcommand{\opU}{\operatorname*{U}}
\newcommand{\refT}[1]{Theorem~\ref{T:#1}}
\newcommand{\refS}[1]{Section~\ref{S:#1}}
\newcommand{\refU}[1]{Subsection~\ref{U:#1}}
\newcommand{\refL}[1]{Lemma~\ref{L:#1}}
\newcommand{\refC}[1]{Chapter~\ref{C:#1}}
\newcommand{\refCL}[1]{Claim~\ref{CL:#1}}
\newcommand{\refP}[1]{Proposition~\ref{P:#1}}
\newcommand{\refD}[1]{Definition~\ref{D:#1}}
\newcommand{\refR}[1]{Remark~\ref{R:#1}}%\renewcommand{\theorem}
%{\arabic{theorem}}%\setcounter{theorem}{1}

%\setcounter{page}{0}

\title
{Regularity for solutions of the two-phase Stefan problem}
\author{Marianne Korten\\
Charles Moore
%\thanks{Math subject classification: }%\thanks{key words:  }
}

\date{}

%\footnote{key words:  }

\maketitle

\begin{abstract}
We consider the two-phase Stefan problem $u_t=\alpha(u)$ where 
$\alpha(u) =u+1$ for $u<-1$, $\alpha(u) =0$ for $-1 \leq u \leq 1$,
and $\alpha(u)=u-1$ for $u \geq 1$.  We show that if $u$ is an $L^2_{loc}$
distributional solution then $\alpha(u)$ is continuous.
\end{abstract}

\section{Introduction}\label{S:1.}
In this paper we will discuss the regularity of weak solutions to the two-phase
Stefan problem
\be\label{e1.1}
\frac{\partial u}{\partial t} =\Delta \alpha(u)
\ee
in a domain $\Omega \subset \R^n \times (0,T)$, for some $T>0$.
Here $\alpha(u)=0$ if $-1 \leq u \leq 1$, $\alpha(u)=u-1$ for $u>1$,
and $\alpha(u)=u+1$ for $u<-1.$

We will show that if $u \in L^2_{loc}(\Omega)$ is a solution in the sense of 
distributions of \eqref{e1.1} (defined precisely below) then $\alpha(u)$ is continuous.
In the case when $u \geq 0$, $u$ is a solution of the one-phase Stefan
problem and Andreucci and Korten \cite{AK} (see also Korten \cite{K}) have shown that
if $u \in L^1_{loc},$ then $\alpha(u)$ is continuous. Although we believe this to be true in the 
two-phase case, we have not been able to obtain this generality and must assume $u \in L^2_{loc}.$

Under the assumption that $u$ is bounded and $\nabla\alpha(u) \in L^2$, Caffarelli and Evans \cite{CE} 
showed that $\alpha(u)$ is continuous. Similar results for more general singular parabolic 
equations were shown by Sacks \cite{S}, Ziemer \cite{Z} and by DiBenedetto \cite{DiB}. 
We will assume these results. We will show that a locally $L^2$ weak solution of \eqref{e1.1}
satisfies the hypotheses of these results (of any of these authors)
to conclude the continuity of $\alpha(u)$. 

Related to this equation is the porous medium equation $u_t=\Delta u^m$, $m>1$.  
This has been studied extensively by many authors, but we mention in particular
the regularity result of Dahlberg and Kenig \cite{DK} who showed that a nonnegative 
$L^m_{loc}$ solution to the porous medium equation is a.e.~ equal to a continuous function.  
The methods in this present paper are descendants (via the work of Andreucci and Korten) 
of the methods of Dahlberg and Kenig found in \cite{DK}.  However, the fact that we are working with 
solutions which can be both positive and negative complicates matters.  To achieve our results we will
perform numerous integrations by parts and cannot determine the sign of the resulting boundary terms
as in the one-phase case.  Consequently we devise a different strategy and introduce new ideas and techniques.

Equation \eqref{e1.1} is a formulation of the two-phase Stefan problem, describing the
flow of heat within a substance which can be in a liquid phase or a solid phase, and for which
there is a latent heat to initiate phase change. This allows for the presence of a ``mushy zone", that is, a 
region which is between the liquid and solid phases.
In this model $u$ represents the enthalpy
and $\alpha(u)$ the temperature.  We have assumed that the thermal conductivity in both the solid and liquid
phases is the same. These conductivities are determined by the slope of the function $\alpha(u)$
in the regions $u \geq 1$, and $u \leq -1$.  The results below all continue to hold (with minor modifications)
if the slope of $\alpha(u)$ differs in these regions. 

We now state our main result.  Suppose $u \in L^2_{loc}(\Omega)$ where $\Omega$ is a domain contained in $\mathbb{R}^n\times (0,T)$.  We consider distributional solutions of the equation $u_t =\Delta \alpha(u)$, that is,
$u$ which satisfy
\bee
\iint_{\Omega} \alpha(u) \Delta \varphi + u \varphi_t dx \, dt =0
\eee
for every $\varphi \in C^{\infty}$ with compact support in $\Omega.$

\begin{theorem}  Suppose $u \in L^2_{loc}(\Omega)$ is a solution of
$u_t=\Delta \alpha(u).$  Then $\alpha(u)$ is a.e.~ equal to a continuous function.
\end{theorem}

We do not expect, in general, such a result for $u$.  As noted in Korten \cite{K1},
the solution to the Cauchy problem $u_t=\Delta \alpha(u)$ on $\mathbb{R}^{n+1}_+$
with initial data $0 \leq u_I(x) \leq 1$ is just $u(x,t)= u_I(x)$. Thus, we cannot expect
$u(x,t)$ to be any smoother than $u_I(x).$

The paper is structured as follows.  In section 2 we prove energy estimates for
weak solutions of the two phase problem.  These show that $\nabla \alpha(u)$
and $\alpha(u)_t$ exist locally in $L^2.$  In section 3, we show that $|\alpha(u)|$
is subcaloric.  An immediate consequence is that $\alpha(u)$ is locally bounded.  This,
combined with the energy estimates and previously mentioned theorem of DiBenedetto \cite{DiB}
(or others mentioned above) gives the continuity of $\alpha(u)$.

Throughout, the letter $C$ will denote a constant which may vary from line to line.

The work of the first author was partially supported by a Kansas EPSCoR grant under 
agreement NSF32169/KAN32170 and a Kansas State University mentoring grant.
The second author would like to thank the National University of Ireland at Galway
for their hospitality during part of this work.

\section{Energy Estimates}\label{S:2.}
We establish that $\alpha(u)$ has derivatives which are locally in $L^2$.

\begin{theorem}\label{T:2.1}
Suppose $\Omega \subseteq \mathbb{R}_+^{n+1}$ and $u\in L^2_\loc(\Omega)$ is a distributional
solution of $u_t=\Delta\alpha(u)$ on $\Omega$. Suppose $r<R$, $T_0<t_0<t_1<T_1$, set
$\omega=(t_0,t_1)\times B(x_0,r)$ and
$\tilde\omega=(T_0,T_1)\times B(x_0,R)$, and suppose the closure of $\tilde\omega$ 
is contained in $\Omega$.  
Then $\nabla\alpha(u)$, $\alpha(u)_t$ exist in $L^2(\omega)$ and there exists a constant $C$, 
depending only on $\omega$ and $\tilde\omega$ such that
\be\label{e2.1}
\iint_{\omega} |\nabla\alpha(u)|^2 dx \, dt\leq C\iint_{\tilde\omega} u^2 dx \, dt\ee
and
\be\label{e2.2}
\iint_{\omega} \left|\frac{\partial}{\partial t}\alpha(u)\right|^2 dx \, dt 
\leq C\iint_{\tilde\omega} u^2 dx \, dt\ee
\end{theorem}

\begin{proof}

Let $\varphi_m(y,s)=\rho_m(y)\tau_m(s)$, $m=1,2,\dots$ where $\rho_m$, $\tau_m$ are 
smooth mollifiers, radial, centered at $0$, compactly supported, and tending to $\delta_0$.
For $(x,t)\in\Omega$ 
and $m$ sufficiently large (depending on $(x,t)$), 
$\varphi_m(x-y,t-s)$ is a test function supported in $\Omega$ and thus
\bee \iint_\Omega u(y,s)\frac{\partial \varphi_m}{\partial t}(x-y,t-s)+\alpha(u(y,s))\Delta \varphi_m(x-y,t-s) dy \, ds=0.
\eee
In the course of the proof, we will need to define three nested
domains between $\omega$ and $\tilde\omega$. To simplify notation, set $\omega_1=\omega$, 
$\omega_5=\tilde\omega$ and we will define $\omega_2$, $\omega_3$ and $\omega_4$ with
$\omega_1 \subset \omega_2 \subset \omega_3 \subset \omega_4 \subset \omega_5.$

For $m=1,2,3,\dots$ set
\bee u_m(x,t)=\iint_\Omega u(y,s)\chi_{\omega_5}(y,s)\varphi_m(x-y,t-s)dy \, ds.\eee
Then for all $(x,t)$ and $m$,
$ |u_m(x,t)|\leq M(u\chi_{\omega_5}(x,t))$
where $M$ denotes the Hardy-Littlewood maximal function (in both the variables $(x,t)$).

Choose $a$, $\frac{3T_0+t_0}{4}<a<\frac{T_0+t_0}{2}$ and $b$, $\frac{t_1+T_1}{2}<b< \frac{t_1+3T_1}{4}$ such that
\be\label{e2.3}
\int_{B(x_0,R)} |M(u\chi_{\omega_5})(x,a)|^2 dx 
\leq C\iint_{\omega_5} |M(u\chi_{\omega_5})(x,t)|^2 dx \, dt .
\ee
with a similar inequality for $b$.

In a similar fashion, set $w_m=\alpha(u)\chi_{\omega_5}\ast\varphi_m$.  Define 
$\omega_4=B\left(x_0,\frac{r+3R}{4}\right)\times\left(\frac{3T_0+t_0}{4},\frac{t_1+3T_1}{4}\right)$.
Then on $\omega_4$, 
$\frac{\partial}{\partial t}u_m-\Delta w_m=0$ for all $m$ sufficiently large.
Using cylindrical coordinates we can choose an $r_1$, $\frac{r+R}{2}<r_1<\frac{r+3R}{4}$ so that
\be\label{e2.4}
\int_{\partial B(x_0,r_1)\times(T_0,T_1)} |M(\alpha(u)\chi_{\omega_5})|^2d\sigma
\leq C\iint_{\omega_5}|M(\alpha(u)\chi_{\omega_5})(x,t)|^2dx \, dt.
\ee
Then by \eqref{e2.3}, for all sufficiently large $m$,
\be\label{e2.5}
\begin{aligned}
\int_{B(x_0,r_1)} u_m(x,a)^2dx & \leq \int_{B(x_0,R)} |M(u\chi_{\omega_5})(x,a)|^2dx\\
   &\leq C\iint_{\omega_5}|M(u\chi_{\omega_5})(x,t)|^2dx \, dt\\
   & \leq C\iint_{\omega_5} u^2dx \, dt
   \end{aligned}\ee
Likewise, by \eqref{e2.4}
\be\label{e2.6}
\begin{aligned}
\int_{\partial B(x_0,r_1)\times(a,b)} w^2_m d\sigma 
   & \leq \int_{\partial B(x_0,r_1)\times(T_0,T_1)} |M(\alpha(u)\chi_{\omega_5})|^2d\sigma\\
   & \leq C\iint_{\omega_5} |M(\alpha(u)\chi_{\omega_5})|^2dx \, dt\\
   & \leq C\iint_{\omega_5} |\alpha(u)|^2dx \, dt.
   \end{aligned}\ee
Let $\alpha_m(s)$ be a smooth regularization of $\alpha(s)$ such that 
$\alpha_m(s)=\alpha(s)$ for $|s|\geq 1+\frac{1}{m}$, 
$\alpha_m(s)$ is strictly increasing and $\alpha_m(s)\not=0$ except for $s=0$.
Put $\omega_3=B(x_0,r_1) \times (a,b).$
Let $v_m$ be a solution to
\bee
\begin{cases}v_t=\Delta\alpha_m(v) &\hbox{\ on\ }\omega_3\\
  v(x,a)=u_m(x,a)               &x\in B(x_0,r_1)\\
 \alpha_m(v)=w_m          &\hbox{\ on\ } \partial B(x_0,r_1)\times(a,b).
  \end{cases}\eee   
Choose $\phi(x)$ so that $\phi=0$ on $\partial B(x_0,r_1)$, $ \Delta\phi=1$ on $B(x_0,r_1)$.
Then $\phi<0$ on $B(x_0,r_1)$ and $\frac{\partial\phi}{\partial n}=c_1>0$ on $\partial B(x_0,r_1)$, 
where $n$ is the outward normal and $c_1$ is a constant depending only on $r_1$ and the dimension.

By Green's theorem we have
\bee\begin{aligned}
&\int_{B(x_0,r_1)} (\alpha_m(v_m))^2\Delta\phi dx \\
  &= \int_{B(x_0,r_1)} \Delta(\alpha_m(v_m))^2\phi dx 
      + \int_{\partial B(x_0,r_1)} (\alpha_m(u_m))^2\frac{\partial\phi}{\partial n}d\sigma
      -\int_{\partial B(x_0,r_1)}\phi\frac{\partial}{\partial n} [\alpha_m(v_m)]^2 d\sigma\\
   &=2\int_{B(x_0,r_1)} \Delta\alpha_m(v_m)\alpha_m(v_m)\phi dx
      +2\int_{B(x_0,r_1)}|\nabla\alpha_m(v_m)|^2\phi dx
       +c_1\int_{\partial B(x_0,r_1)}\alpha_m(v_m)^2d\sigma\\
    & \leq 2\int_{B(x_0,r_1)} v_{m_t}\alpha_m(v_m)\phi dx
       +c_1\int_{\partial B(x_0,r_1)} \alpha_m(v_m)^2 d\sigma\\
     &=2\frac{d}{dt}\int_{B(x_0,r_1)} A_m(v_m)\phi dx
        +c_1\int_{\partial B(x_0,r_1)} \alpha_m(v_m)^2d\sigma
       \end{aligned}\eee
where $A_m$ is an antiderivative of $\alpha_m$. Integrate from $a$ to $b$ to obtain
\bee\begin{aligned}
\int^b_a\int_{B(x_0,r_1)}
   \alpha_m(v_m)^2 \Delta\phi dx \, dt\leq  2 \int_{B(x_0,r_1)}& A_m(v_m(x,b))\phi dx 
    -2\int_{B(x_0,r_1)} A_m(v_m(x,a))\phi dx \\ 
    & +\int^b_a\int_{\partial B(x_0,r_1)} \alpha_m(v_m)^2d\sigma dt.
     \end{aligned}\eee
Now $0\leq A_m(x)\leq x^2$, $\phi<0$, $\Delta\phi=1$ and recalling $\omega_3=B(x_0,r_1)\times(a,b)$, this
yields:
\be\label{e2.7}\begin{aligned}
  \iint_{w_3}\alpha_m(v_m)^2 dx \, dt
    & \leq C\int_{B(x_0,r_1)}v_m(x,a)^2dx 
      +\int^b_a\int_{\partial B(x_0,r_1)} \alpha_m(v_m)^2d\sigma dt \\
     &=C\int_{B(x_0,r_1)} u_m(x,a)^2dx+\int^b_a\int_{\partial B(x_0,r_1)} w^2_m d\sigma dt\\
      & \leq C\iint_{w_5}u^2dx \, dt
     \end{aligned}\ee
where for the last inequality we have used \eqref{e2.5} and \eqref{e2.6}.
Let $\psi(x)$ be a nonnegative $C^\infty_0(\R^n)$ function such that
$\psi\equiv 1$ on $B(x_0,\frac{3r+R}{4})$, $\psi\equiv 0$ outside $B(x_0,\frac{r+R}{2})$. 
To simplify notation set $B(x_0,\frac{3r+R}{4})=B_1$, $B(x_0,\frac{r+R}{2})=B_2$.
Then     
\bee\begin{aligned}
\int_{B_2}\psi\alpha_m(v_m){v_m}_tdx
  &=\int_{B_2}\psi\alpha_m(v_m)\Delta\alpha_m(v_m)dx\\
  &=-\int_{B_2}\nabla\psi\cdot\nabla\alpha_m(v_m) \alpha_m(v_m)dx
      -\int_{B_2}\psi |\nabla\alpha_m(v_m)|^2dx\\
    &=\frac{1}{2}\int_{B_2}\Delta\psi\alpha_m(v_m)^2dx 
     - \int_{B_2}\psi|\nabla\alpha_m(v_m)|^2dx.
   \end{aligned}\eee
Rearrange and integrate from $a$ to $b$ to obtain
 \be\label{e2.8}
 \begin{aligned}
 \int^b_a\int_{B_2}\psi (\nabla\alpha_m(v_m))^2 dx\, dt
   & =\frac{1}{2}\int^b_a\int_{B_2}\Delta\psi |\alpha_m(v_m)^2| dx \, dt
     -\int^b_a\int_{B_2}\psi\frac{d}{dt}A_m(v_m)dx \, dt\\
   &\leq C\int^b_a\int_{B_2}\alpha_m(v_m)^2dx \, dt +\int_{B_2}\psi(x)A_m(v_m(x,a))dx\\
   &\leq C\int^b_a\int_{B_2}\alpha_m(v_m)^2dx \, dt + \int_{B_2}v_m(x,a)^2dx\\
   &\leq C\iint_{w_3} \alpha_m(v_m)^2dx \, dt +\int_{B(x_0,r_1)}v_m(x,a)^2dx\\
   &\leq C\iint_{w_5}u^2dx \, dt 
   \end{aligned}\ee
where we have used \eqref{e2.7} and \eqref{e2.5} and the definition of $v_m$ for the last inequality.

We now seek a similar estimate for the $t$ derivative.
Let $\eta(x)$ be a nonnegative $C^\infty_0(\R^n)$ function such that $\eta\equiv 1$ on $B(x_0,r)$, $\eta \equiv 0$ 
outside $B_1$ and so that $\|\frac{\nabla \eta}{\sqrt{\eta}}\|_\infty<\infty$. Note that
$\alpha_m(v_m)_t=\alpha'_m (v_m) {v_m}_t$ and $0< \alpha_m' \leq 1$ so that
$(\alpha_m(v_m)_t)^2 \leq \alpha_m(v_m)_t v_{m_t}.$  Then
\bee\begin{aligned}
\int_{B_1}\eta (\alpha_m(v_m)_t)^2dx
  &\leq \int_{B_1}\eta \alpha_m(v_m)_t v_{m_t}dx\\
  &=\int_{B_1} \eta \alpha_m(v_m)_t\Delta\alpha_m(v_m)dx\\
  &=-\int_{B_1}\nabla \eta\cdot\nabla\alpha_m(v_m)\alpha_m(v_m)_t dx
      -\int_{B_1}\eta\nabla\alpha_m(v_m)_t\cdot\nabla\alpha_m(v_m)dx\\
    &=-\int_{B_1}\nabla \eta\cdot\nabla\alpha_m(v_m)\alpha_m(v_m)_tdx
        -\frac{1}{2}\frac{d}{dt}\int_{B_1} \eta|\nabla\alpha_m(v_m)|^2dx
  \end{aligned}\eee
Integrate from $c$ to $d$, where $c$ and $d$ are to be chosen momentarily.
We obtain
\be\label{e2.9}
\begin{aligned}
\int^d_c \int_{B_1} & \eta(\alpha_m(v_m)_t)^2dx \, dt\\
  &\leq \int^d_c\left|\int_{B_1} \sqrt{\eta}\frac{\nabla \eta}{\sqrt{\eta}}
     \nabla\alpha_m(v_m)\alpha_m(v_m)_tdx\right| dt
     +\frac{1}{2}\int_{B_1}\eta(x)|\nabla\alpha_m(v_m)(x,c)|^2dx\\
   &\leq\left\|\frac{\nabla \eta}{\sqrt{\eta}}\right\|_\infty 
      \left( \int^d_c\int_{B_1}|\nabla\alpha_m(v_m)|^2dx \, dt\right)^{\frac{1}{2}}
      \left(\int^d_c\int_{B_1}\eta(\alpha_m(v_m)_t)^2dx \, dt\right)^{\frac{1}{2}}\\
    &\qquad\qquad +\frac{1}{2}\int_{B_1}\eta(x)|\nabla\alpha_m(v_m)(x,c)|^2dx.
    \end{aligned}\ee
 Choose $c_m$ (depending on $m$),  
$\frac{T_0+3t_0}{4}<c_m<t_0$, so that
\be\label{e2.10}
\int_{B_1}\eta(x)|\nabla\alpha_m(v_m)(x,c_m)|^2dx  \leq C \int_a^b\int_{B_2}\psi|\nabla\alpha_m(v_m)|^2dx \, dt.
\ee
Put $d=t_1$, $c=c_m$ in \eqref{e2.9}. Then recalling that $\psi \equiv 1$ on $B_1$, and using \eqref{e2.10} and \eqref{e2.8} we have
\bee \begin{aligned}
\int^{t_1}_{c_m}&\int_{B_1}  \eta(\alpha_m(v_m)_t)^2dx \, dt \\
   &\hspace{-.25in}
   \leq \left\|\frac{\nabla \eta}{\sqrt{\eta}}\right\|_\infty
    \left(\int^{t_1}_{c_m}\int_{B_2}\psi |\nabla\alpha_m(v_m)|^2  dx\,dt \right)^{\frac{1}{2}}
      \cdot \left(  \int^{t_1}_{c_m}\int_{B_1}\eta(\alpha_m(v_m)_t)^2dx \, dt \right)^{\frac{1}{2}}
  + C\iint_{\omega_5}u^2dx\\
  & \hspace{-.25in}
    \leq\left\|\frac{\nabla \eta}{\sqrt{\eta}}\right\|_\infty
    \left( \iint_{\omega_5} u^2dx \, dt\right)^{\frac{1}{2}}
     \left( \int^{t_1}_{c_m}\int_{B_1}\eta(\alpha_m(v_m)_t)^2dx \, dt \right)^{\frac{1}{2}}
      +C\iint_{\omega_5}u^2dx \, dt
\end{aligned}\eee 
from which it follows that
 \bee \int^{t_1}_{c_m}\int_{B_1}\eta(\alpha_m(v_m)_t)^2dx \, dt \leq C \iint_{\omega_5}u^2dx \, dt \eee
 and consequently
 \be\label{e2.11}
 \iint_{\omega_1}(\alpha_m(v_m)_t)^2dx \, dt \leq C\iint_{\omega_5}u^2dx \, dt .
 \ee
 Thus, recalling $\omega=\omega_1$, $\tilde\omega=\omega_5$,  \eqref{e2.8} and \eqref{e2.11} give 
 \be\label{e2.12}
  \iint_\omega|\nabla\alpha_m(v_m)|^2dx \, dt\leq C\iint_{\tilde\omega} u^2dx \, dt \ \   \text{and} \ 
  \iint_\omega (\alpha_m(v_m)_t)^2dx \, dt \leq C\iint_{\tilde\omega} u^2 dx \, dt .
 \ee
To obtain \eqref{e2.1} and \eqref{e2.2} we will need to take limits.  We first remark that
with more care, similar estimates could be obtained with any compact set
$K\subset\omega_3$ replacing $\omega=\omega_1$ on the left hand side of the inequalities in
\eqref{e2.12}; naturally, the constants on the right hand side depend on the position of $K$ within $\omega_3.$
Thus, from \eqref{e2.7} and this observation, we have:
\be\label{e2.13}\begin{aligned}
\iint_{w_3}\alpha_m(v_m)^2dx \, dt \leq C\iint_{\tilde\omega}u^2dx \, dt ,& \qquad \iint_K |\nabla\alpha_m(v_m)|^2 dx \, dt \leq C(K)\iint_{\tilde\omega} u^2 dx \, dt\\
\text{and } \iint_K \left|\frac{\partial}{\partial t}\alpha_m(v_m)\right|^2 & dx \, dt \leq C(K)\iint_{\tilde\omega} u^2 dx \, dt
\end{aligned}\ee
for every compact $K \subset \omega_3.$

By Rellich-Kondrachov there exists a subsequence $\{\alpha_{m_k}(v_{m_k})\}$ of $\{\alpha_m(v_m)\}$
(which we still write as $\{\alpha_m(v_m)\}$) and $h\in L^2(\omega_3)$ such that 
$\alpha_m(v_m) \to h$ in $L^2(K)$ for every compact set $K \subset \omega_3.$
By taking subsequences, if necessary, we also may assume this convergence is a.e.
By weak compactness, and again, by taking subsequences, we may assume that
$\alpha_m(v_m) \to h$ weakly in $L^2(\omega_3)$. Equation \eqref{e2.13} implies that the
$L^2(\omega_3)$ norms of the $v_m$ are uniformly bounded, hence there exists a subsequence,
(still denoted by $v_m$) such that $v_m \to v \in L^2(\omega_3)$ weakly.

We claim that $\alpha(v)=h$.  First note that $\|\alpha_m -\alpha\|_{\infty} \to 0$,
so that for a.e $x \in \omega_3$, $\alpha(v_m) \to h$.  Consider the set where $h>0$.
Then for a.e $x$ in this set, $\alpha(v_m(x)) \to h(x) >0$, and hence $v_m(x) \to h(x)+1$.
Thus, $v(x)=h(x)+1$ for a.e $x$ in the set where $h(x)>0$.
Similarly, on $h<-1$, $v(x)= h(x)-1$ a.e.  On the set $h(x)=0$ we must have 
$-1 \leq \liminf v_m(x) \leq \limsup v_m(x) \leq 1$ a.e. To see this consider an $x$ at which there exists
a subsequence $v_{m_k}(x)$ which converges to $y_0 \notin [-1,1].$ Then for this $x$,
$\alpha(v_{m_k}(x)) \to \alpha(y_0) \neq 0$ which implies $\alpha(v_m(x))\nrightarrow h(x)$.
Thus, $-1 \leq \liminf v_m(x) \leq \limsup v_m(x) \leq 1$ a.e. on $h(x)=0$, and hence
$-1 \leq v(x) \leq 1$ a.e. on $h=0$.  We conclude that $\alpha(v)=h$ a.e.

Summarizing, we have $v_m \to v$ weakly in $L^2(\omega_3)$ and $\alpha_m(v_m) \to \alpha(v)$
weakly in $L^2(\omega_3)$, a.e. on $\omega_3$ and in $L^2(K)$ for every compact subset 
$K$ of $\omega_3$. To finish the proof we show that $\alpha(u)=\alpha(v)$ a.e. on $\omega_3$. 
Using integration by parts, and recalling that $\alpha_m(v_m)=w_m$ on $\partial B(x_0,r_1) \times (a,b)$, we compute
\be \label{e2.14}
\begin{aligned}
\iint_{\omega_3}  (v_m - u_m) & (\alpha_m(v_m) -w_m) dx \, dt  \\
& =\int_a^b \int_{B(x_0,r_1)} \int_a^t {v_m}_t(x, \tau) - {u_m}_t(x,\tau)d\tau (\alpha_m(v_m(x,t))-w_m(x,t)) dx \, dt   \\
&=-\int_a^b\int_{B(x_0,r_1)}  \int_a^t \nabla(\alpha_m(v_m) -w_m) (x,\tau) d\tau \nabla (\alpha_m(v_m)-w_m)dx \, dt \\
&=-\frac{1}{2} \int_a^b \int_{B(x_0, r_1)} \frac{d}{dt} \left|\int_a^t \nabla(\alpha_m(v_m) -w_m) d\tau\right|^2 dx \, dt \\
&=-\frac{1}{2} \int_{B(x_0,r_1)} \left|\int_a^b \nabla(\alpha_m(v_m)-w_m) d\tau\right|^2 dx \leq 0
\end{aligned}
\ee
We need to take limits as $m \to \infty$ in this inequality. Write
\bee
\begin{aligned}
&\iint_{\omega_3}  (v_m - u_m)(\alpha_m(v_m) -w_m) dx \, dt = \iint_{\omega_3}v_m \alpha_m(v_m) dx\, dt 
 + \iint_{\omega_3}v_m (-w_m)dx \, dt \\ 
& \ \qquad + \iint_{\omega_3}(-u_m)(\alpha_m(v_m))dx\, dt + \iint_{\omega_3}u_m w_m dx\, dt= I + II + III + IV
\end{aligned}
\eee
Since $u_m \to u$ a.e. and in $L^2(\omega_3)$, $w_m \to \alpha(u)$ a.e. and in $L^2(\omega_3)$,
$v_m \to v$ weakly, and $\alpha_m(v_m) \to \alpha(v)$ weakly, we conclude 
\bee \begin{aligned}
II \to \iint_{\omega_3} v(-\alpha(u))&  dx \, dt, \qquad III \to \iint_{\omega_3} (-u)(\alpha(v))dx \, dt, \\
& \text{and } \  IV \to \iint_{\omega_3} u \alpha(u) dx \, dt.
\end{aligned}
\eee
Expand out
\bee
\iint_{\omega_3} (\alpha_m(v_m)-\alpha_m(v))(v_m-v)dx \geq 0,
\eee
take $m \to \infty$ (make use of the fact that $ \|\alpha_m -\alpha\|_{\infty} \to 0$) to conclude
\bee
\liminf_{m \to \infty} \iint_{\omega_3} \alpha_m(v_m)v_m dx \geq \iint_{\omega_3} \alpha(v)v dx.
\eee
This combined with the estimates for II-IV and \eqref{e2.14} yields
\bee
\iint_{\omega_3} (v-u)(\alpha(v)-\alpha(u)) dx \, dt \leq 0.
\eee
Since the integrand of this is nonnegative, we conclude $\alpha(u)=\alpha(v)$ a.e.~ on
$\omega_3$.  This completes the proof of the theorem.
\end{proof}

\section{$|\alpha(u)|$ is subcaloric}

\begin{theorem}\label{T:3.1}  $|\alpha(u)|$ is weakly subcaloric, that is, it satisfies
\bee
 \int_{\Omega}-\nabla|\alpha(u)|\nabla\eta + |\alpha(u)|\eta_t dx \, dt\geq 0
\eee
for any nonnegative $\eta \in W^{1,2}_0 (\Omega).$

\end{theorem}

\begin{proof}
Let $0\leq \eta  \in W^{1,2}_0 (\Omega).$
For $h>0$ set
\bee
\phi_h(x)=
\begin{cases} 1 &\hbox{if\ } x>h\\
\frac{2}{h}x-1 & \frac{h}{2} \leq x<h\\
  0&\hbox{if\ } |x| < \frac{h}{2}\\
  \frac{2}{h}x+1 & -h<x\leq -\frac{h}{2}\\
 -1 & \hbox{if\ } x<-h.
  \end{cases}
\eee
Then $\eta\phi_h(\alpha(u))$ is supported in
$\left\{|\alpha(u)|> \frac{h}{2}\right\}$ and thus,
$\iint (\alpha(u)-u) [\eta\phi_h(\alpha(u))]_t dx \, dt =0.$

Then
\be\label{e3.1}
\begin{aligned} 
0&=\iint \alpha(u)[\eta\phi_h(\alpha(u))]_t - \nabla\alpha(u)
   \nabla[\eta\phi_h(\alpha(u))]dx\,dt  \\
 &=\iint\alpha(u) \eta_t\phi_h(\alpha(u)) dx\,dt
   + \iint\alpha(u)\eta\phi'_h
   (\alpha(u)) \alpha(u)_t  dx\,dt   \\
 &\quad - \iint\nabla\alpha(u) \nabla \eta\phi_h(\alpha(u)) dx\,dt
   -\iint\nabla\alpha(u) \eta\phi'_h(\alpha(u)) \nabla\alpha(u) dx\,dt \\
    & =I+II+III+IV  
\end{aligned}
\ee

We investigate each of these as $h\to 0$.
As $h\to 0$, $\phi_h(\alpha(u))\to \hbox{sgn}(\alpha(u))$ so that
$I\to\iint|\alpha(u)|n_tdx\,dt$.
To estimate $II$, first note that
$$\phi'_h(\alpha(u))=\frac{2}{h} \chi_{\left\{\frac{h}{2} < |\alpha|< h\right\}}.$$
Then
\bee
 II=\iint\alpha(u)\eta\frac{2}{h} \chi_{\{\frac{h}{2} <  |\alpha(u)|<h\}}
  \frac{\partial}{\partial t} \alpha(u) dx\,dt 
\eee
and consequently, 
\bee
  |II|\leq \iint|\alpha(u)|\eta\frac{2}{h} \chi_{\{\frac{h}{2}
    < |\alpha(u)|<h\}}
   \left|\frac{\partial}{\partial t}\alpha(u)\right| dx\,dt
   \leq 2\iint \eta \chi_{\{\frac{h}{2} < \alpha(u)<h\}}
    \left |\frac{\partial}{\partial t}(\alpha(u))\right| dx\,dt.
\eee
Since $\frac{\partial}{\partial t}\alpha(u)
\in L^2_{loc}(\Omega)$, $II\to 0$ as $h\to 0$

To estimate $III$, note that when $|\alpha(u)|>h$,
$\nabla\alpha(u) \phi_h (\alpha(u))=\nabla|\alpha(u)|$.
And when $\frac{h}{2} < |\alpha(u)|<h$,
\bee |\nabla \alpha(u) \phi_h(\alpha(u))| \leq |\nabla\alpha(u)|
  \left(\frac{2}{h}|\alpha(u)| + 1\right)
  \eee
Consequently,
\bee\begin{aligned}
 \left| \iint_{\{\frac{h}{2} < |\alpha(u)|<h\}}
   \nabla\alpha(u)\nabla \eta\phi_h(\alpha(u)) dx\,dt\right| & \leq  \iint_{\{\frac{h}{2} < |\alpha(u)|<h\}}
   \left|\nabla\alpha(u)\right|
   \left[ \frac{2}{h} |\alpha(u)|+ 1\right]
   |\nabla \eta| dx\,dt\\
  &\quad \to 0 \hbox{\ as\ } h\to 0 \hbox{\ since\ } \nabla\alpha(u)\in 
L^2_{loc}(\Omega)).
  \end{aligned} 
 \eee
Therefore, as $h\to 0$,
$III\to -\iint \nabla|\alpha(u)|\nabla\eta dx\,dt.$
Note that we can write $IV$ as
\bee 
IV =\iint|\nabla\alpha(u)|^2\eta\ \phi'_h(\alpha(u))dx\,dt
   =\iint|\nabla\alpha(u)|^2\eta\frac{2}{h} \chi_{\left\{\frac{h}{2}<|\alpha(u)|<h\right\}} dx\,dt 
\eee
Thus, letting $h\to 0$ in \eqref{e3.1} yields:
\bee 
0 =\iint|\alpha(u)|\eta_t dx\,dt
      -\iint \nabla|\alpha(u)|\nabla\eta dx\,dt
    -\lim_{h\to 0} \iint|\nabla\alpha(u)|^2
     \frac{2}{h} \chi_{\left\{\frac{h}{2} < |\alpha(u)|<h\right\}}\eta dx\,dt 
\eee
from which the theorem follows.

\end{proof}
\noindent{\bf Remark.} Suppose that instead of $\phi_h$ as defined above, we defined
$$ \phi_h(x)=
  \begin{cases}   1 &\hbox{if\ } x>h\\
  \frac{2}{h}x-1 &\hbox{if\ } \frac{h}{2} <x<h\\
  0& \hbox{otherwise}.
  \end{cases} $$
Then following the computations as in \eqref{e3.1} we 
obtain \eqref{e3.1} with this version of $\phi_h$.
In this case $I\to\iint\alpha(u)^+\eta_t dx\,dt$ and 
as in the above case, $II\to 0$, 
and $III \to  -\iint\nabla\alpha(u)^+\nabla \eta dx\,dt.\\$
Now we may write
$$ IV=\iint |\nabla\alpha(u)^+|^2
 \eta\frac{2}{h} \chi_{\left\{\frac{h}{2}<\alpha(u)<h\right\}} dx\,dt.$$
We obtain:
\bee
 0= \iint\alpha(u)^+\eta_tdx\,dt
    -\iint\nabla\alpha(u)^+\nabla \eta dx\,dt
  -\lim_{h\to 0} \iint|\nabla\alpha(u)^+|^2
    \frac{2}{h} \chi_{\left\{\frac{h}{2} <\alpha(u)<h\right\}}\eta dx\,dt
\eee
Thus, $\alpha(u)^+$ is subcaloric.
In a similar fashion, we may use the function
$$\phi_h(x)= \begin{cases}
 \frac{2}{h}x +1 & -h<x<-\frac{h}{2}\\
 -1 & \hbox{if\ } x<-h\\
 0& \hbox{otherwise}
 \end{cases} $$
and computations such as those above to obtain
\bee 
  0=\iint\alpha(u)^-\eta_t dx\,dt
     -\iint\nabla\alpha(u)^- \nabla \eta dx\,dt
    + \lim_{h\to 0} \iint |\nabla\alpha(u)^-|^2 \frac{2}{h}
    \chi_{\left\{-h<\alpha(u) <-\frac{h}{2}\right\}}\eta dx\,dt
\eee
to conclude that $\alpha(u)^-$ is supercaloric.


\begin{thebibliography}{1000} %number of characters of longest bibitem label

\bibitem[AnKo]{AK} D.~ Andreucci and M.~ K.~ Korten, {\it Initial traces of solutions
to a one-phase Stefan problem in an infinite strip,} Rev.~ Mat.~ 
Ibero\-a\-me\-ri\-cana  9 (2) (1993), 315-332.

\bibitem[CaE]{CE} L.~ A.~ Caffarelli and L.~ C.~ Evans, {\it Continuity of the temperature in the two-phase Stefan problem,}  Arch. Rational Mech. Anal. 81 (1983), 199-220.
 
\bibitem[DiB]{DiB} E.~ DiBenedetto, {\it Continuity of weak solutions to certain singular parabolic equations,} 
Ann. Mat. Pura Appl. (IV) 130 (1982), 131--176. 

\bibitem[DK]{DK}  B.~ E.~ J.~ Dahlberg and C.~ E.~ Kenig, {\it Weak solutions of the porous medium equation,}
Trans. Amer. Math. Soc., 336 (2) (1993), 711-725.

\bibitem[Ko]{K} M.~ K.~ Korten, {\it Non-negative solutions of $u_t = \Delta
(u-1)_+$: Regularity and uniqueness for the Cauchy problem,} Nonl.~ Anal., Th., Meth.~ and Appl. 27 (5) (1996), 589-603. 

\bibitem[Ko1]{K1} M.~ K.~ Korten, {\it A Fatou theorem for the equation $u_t =
\Delta (u-1)_+$,} Proc.~ Am.~ Math.~ Soc. 128 (2) (2000), 439-444.

\bibitem[S]{S} P.~E.~ Sacks, {\it Continuity of solutions of a singular parabolic equation,} 
Nonlinear Anal.~ 7 (4) (1983), 387--409.

\bibitem[Z]{Z} W.~P.~ Ziemer, {\it Interior and boundary continuity of weak solutions of degenerate parabolic equations,} Trans.~ Amer.~ Math.~ Soc. 271 (2) (1982), 733--748.

\end{thebibliography}
\end{document}